\documentclass[12pt]{article}

\newtheorem{defi}{Definition}

\newtheorem{tm}{Theorem}

\newtheorem{rem}{Remark}

\usepackage{latexsym,amssymb}
\usepackage{float}
\usepackage{xcolor}
\usepackage{geometry} 

\begin{document}

\vspace*{0.5cm}

\begin{center}
{\Large\bf  On some distance-regular graphs with many vertices} 
\end{center}

\vspace*{0.5cm}

\begin{center}
        Dean Crnkovi\'c (deanc@math.uniri.hr)\\
		Sanja Rukavina (sanjar@math.uniri.hr)\\
		and \\
		Andrea \v Svob (asvob@math.uniri.hr)\\[3pt]
		{\it\small Department of Mathematics} \\
		{\it\small University of Rijeka} \\
		{\it\small Radmile Matej\v ci\'c 2, 51000 Rijeka, Croatia}\\
\end{center}

\vspace*{0.5cm}

\begin{abstract}
We construct distance-regular graphs, including strongly regular graphs, admitting a transitive action of the Chevalley groups 
$G_2(4)$ and $G_2(5)$, the orthogonal group $O(7,3)$ and the Tits group $T=$$^2F_4(2)'$. Most of the constructed graphs have more than 1000 vertices,
and the number of vertices goes up to 28431. 
Some of the obtained graphs are new. 

\end{abstract}

\bigskip

\noindent {\bf AMS classification numbers:} 05E30, 05E18.


\section{Introduction}

A construction of distance-regular graphs (DRGs), and especially strongly regular graphs (SRGs), from finite groups gave an important contribution to the graph theory and the design theory
(see \cite{crc-srg, drg, low-rank}). 
Recently, in \cite{crsSRG} the authors found new SRGs admitting a transitive action of some finite simple groups. 
Research presented in this paper show how one can use groups as a tool to produce wide range of regular graphs. 
Moreover, in that way an insight to particular groups is given and some DRGs (and especially SRGs) with many vertices are constructed. 
In literature, and also implemented in Sagemath (see \cite{dima, sage}), there are a few examples of SRGs and DRGs  with many vertices that do not belong to some infinite 
family of graphs. 
The expression "many vertices" in this paper refers to the number of vertices that is out of range given in \cite{aeb} and \cite{drg}. 
For the sake of completeness of the presented results we also include obtained graphs on less than 1000 vertices. 

We assume that the reader is familiar with the basic facts of the group theory, the theory of strongly regular graphs and the theory of distance-regular graphs. 
We refer the reader to \cite{atlas, r} for relevant background reading in the group theory, to \cite{crc-srg, tonchev-book} for the theory of strongly regular graphs, 
and to \cite{drg,dam} for the theory of distance-regular graphs. 

In this paper we study the exceptional groups of Lie type $G_2(4)$ and $G_2(5)$, the orthogonal group $O(7,3)$ and the Tits group $T$, which are the simple groups of orders 
251596800, 5859000000, 4585351680 and 17971200, respectively. 
We refer the reader to \cite{atlas, wilson} for more details about these groups.
 
Using the method outlined in Section \ref{SRG_groups} we constructed and classified SRGs and DRGs of diameter $d \ge 3$ from above mentioned simple groups as follows:

\begin{itemize}
	\item up to 20000 vertices admitting a primitive action of the group $G_2(4)$, and up to 13000 vertices admitting an imprimitive action of the group $G_2(4)$,
	\item up to 400000 vertices admitting a primitive action of the group $G_2(5)$, and up to 15000 vertices admitting an imprimitive action of the group $G_2(5)$,
	\item up to 250000 vertices admitting a primitive action of the group $O(7,3)$, and up to 10000 vertices admitting an imprimitive action of the group $O(7,3)$,
	\item up to 12000 vertices admitting a primitive action of the group $T$, and up to 7000 vertices admitting an imprimitive action of the group $T$.
	\end{itemize}

Based on \cite{crc-srg, aeb, LintBrouwer} and the other literature we conclude that the constructed strongly regular graphs with parameters $(28431,3150,621,315)$ and 
$(28431,2880,324,288)$ are the first known examples of SRGs with these parameters.

The classification, that was conducted by use of computers, was computationally demanding. 
The running time complexity of the algorithm used for the construction of the graphs depends on the number
of parameters, such as the size of the used subgroup, the number of orbits of a block stabilizer (in $G_2(4)$, $G_2(5)$, $O(7,3)$ or $T$), 
the number of vertices of the graphs. 

Note that in the paper we use the notation from \cite{atlas} to describe the structures of the groups.
To find the graphs and compute their full automorphism groups, we used programmes written for Magma \cite{magma} and GAP \cite{GAP2016}.
The constructed SRGs and DRGs can be found at the link: 
\begin{verbatim}
http://www.math.uniri.hr/~asvob/DRGs_manyVertices.zip.
\end{verbatim}

\section{Preliminaries}

In this section we define coherent configurations and association schemes, which are the tools for the construction of graphs presented in this paper.
We also give basic definitions and properties of DRGs and SRGs.

\begin{defi}\label{CohConf}
A coherent configuration on a finite non-empty set $\Omega$ is an ordered pair $(\Omega, \mathcal{R})$ with $\mathcal{R}=\{R_0,R_1,\dots,R_d\}$ a set of non-empty relations on $\Omega$, 
such that the following axioms hold.

\begin{itemize}
	\item[(i)] $\displaystyle\sum_{i=0}^t R_i$ is the identity relation, where $\{R_0,R_1,\dots,R_t\} \subseteq \{R_0,R_1,\dots,R_d\}$.
		\item[(ii)] $\mathcal{R}$ is a partition of $\Omega^2$.
		\item[(iii)] For every relation $R_i\in \mathcal{R}$, its converse $R_i^T=\{(y,x):(x,y)\in R_i\}$ is in $\mathcal{R}$.
		\item[(iv)] There are constants $p_{ij}^k$ known as the intersection numbers of the coherent configuration $\mathcal{R}$, such that for $(x,y)\in R_k$, the number of elements $z$ in $\Omega$ 
		for which $(x,z)\in R_i$ and $(z,y)\in R_j$ equals $p_{ij}^k$.
\end{itemize}
\end{defi}

We say that a coherent configuration is homogeneous if it contains the identity relation, i.e., if $R_0=I$.
If $\mathcal{R}$ is a set of symmetric relations on $\Omega$, then a coherent configuration is called symmetric. 
A symmetric coherent configuration is homogeneous (see \cite{cameron}).
Symmetric coherent configurations are introduced by Bose and Shimamoto in \cite{bose-associate} and called association schemes.
An association scheme with relations $\{R_0,R_1,\dots,R_d\}$ is called a $d$-class association scheme.


Let $\Gamma$ be a graph with diameter $d$, and let $\delta(u,v)$ denote the distance between vertices $u$ and $v$ of $\Gamma$.
The $i$th-neighborhood of a vertex $v$ is the set $\Gamma_{i}(v) = \{ w : d(v,w) = i \}$. 
Similarly, we define $\Gamma_{i}$ to be the $i$th-distance graph of $\Gamma$, that is, the vertex set of $\Gamma_{i}$ is the same as for $\Gamma$, with adjacency in $\Gamma_{i}$ defined by the $i$th distance relation in $\Gamma$.
We say that $\Gamma$ is distance-regular if the distance relations of $\Gamma$ give the relations of a $d$-class association scheme, that is, for every choice of $0 \leq i,j,k \leq d$, all vertices $v$ and $w$ with $\delta(v,w)=k$ satisfy $|\Gamma_{i}(v) \cap \Gamma_{j}(w)| = p^{k}_{ij}$ for some constant $p^{k}_{ij}$.

In a distance-regular graph, we have that $p^{k}_{ij}=0$ whenever $i+j < k$ or $k<|i-j|$.
A distance-regular graph $\Gamma$ is necessarily regular with degree $p^{0}_{11}$; more generally, each distance graph $\Gamma_{i}$ is regular with degree $k_{i}=p^{0}_{ii}$.

An equivalent definition of distance-regular graphs is the existence of the constants $b_{i}=p^{i}_{i+1,1}$ and $c_{i}= p^{i}_{i-1,1}$ for $0 \leq i \leq d$ (notice that $b_{d}=c_{0}=0$).
The sequence $\{b_0,b_1,\dots,b_{d-1};c_1,c_2,\dots,c_d\}$, where $d$ is the diameter of $\Gamma$ is called the intersection array of $\Gamma$. 
Clearly, $b_0=k$, $b_d=c_0=0$, $c_1=0$.

A regular graph is strongly regular of type $(v,k, \lambda , \mu )$ if it has $v$ vertices, degree $k$,
and if any two adjacent vertices are together adjacent to $\lambda$ vertices, while any two non-adjacent vertices are together adjacent to $\mu$ vertices.
A strongly regular graph  of type $(v,k, \lambda , \mu )$ is usually denoted by SRG$(v, k, \lambda, \mu)$. 
A strongly regular graph is a distance-regular graph with diameter 2 whenever $\mu \neq 0$.
The intersection array of an SRG is given by $\{k,k-1-\lambda;1,\mu\}$.

\section{SRGs and DRGs constructed from the groups} \label{SRG_groups}

In this section we give a classification of SRGs and DRGs of diameter $d \ge 3$ as follows:
\begin{itemize}
	\item up to 20000 vertices admitting a primitive action of the group $G_2(4)$, and up to 13000 vertices admitting an imprimitive action of the group $G_2(4)$,
	\item up to 400000 vertices admitting a primitive action of the group $G_2(5)$, and up to 15000 vertices admitting an imprimitive action of the group $G_2(5)$,
	\item up to 250000 vertices admitting a primitive action of the group $O(7,3)$, and up to 10000 vertices admitting an imprimitive action of the group $O(7,3)$,
	\item up to 12000 vertices admitting a primitive action of the group $T$, and up to 7000 vertices admitting an imprimitive action of the group $T$.
\end{itemize}

\bigskip

Let $G$ be a finite permutation group acting on the finite set $\Omega$. This action induce the action of the group $G$ on the set $\Omega\times \Omega$. For more information see \cite{wielandt}. 
The orbits of this action are the sets of the form $\{(\alpha g,\beta g): g\in G\}$. If $G$ is transitive, then $\{(\alpha ,\alpha):\alpha\in \Omega\}$ is one such orbit. 
If the rank of $G$ is $r$, then it has $r$ orbits on $\Omega\times\Omega$.
Let $|\Omega|=n$ and $\Delta_i$ is one of these orbits. We say that the $n\times n$ matrix $A_i$, with rows and columns indexed by $\Omega$ and entries 

$$
A_i(\alpha,\beta)= \left \lbrace
\begin{array}{ll}
1,      & \mathrm{if} \ (\alpha,\beta)\in \Delta_i\\
0, & \mathrm{otherwise}.
\end{array}  \right.
$$
is called the adjacency matrix for the orbit $\Delta_i$.

Let $A_0,\dots,A_{r-1}$ be the adjacency matrices for the orbits of $G$ on $\Omega\times\Omega$. These satisfy the following conditions.

\begin{itemize}
\item[(i)] $A_0=I$, if $G$ is transitive on $\Omega$. If $G$ has $s$ orbits on $\Omega$, then $I$ is a sum of $s$ adjacency matrices.
\item[(ii)] $\displaystyle\sum_{i}A_i=J$, where $J$ is the all-one matrix.
\item[(iii)] If $A_i$ is an adjacency matrix, then so is its transpose $A_i^T$.
\item[(iv)] If $A_i$ and $A_j$ are adjacency matrices, then their product is an integer-linear-combination of adjacency matrices.
\end{itemize}

If $A_i$ is symmetric, then the corresponding orbit is called self-paired. Further, if $A_i=A_j^T$, then the corresponding orbits are called mutually paired.

The graphs obtained in this paper are constructed using the method described in \cite{cms} which can be rewritten in terms of coherent configurations in the following way.

\begin{tm}\label{main}
Let $G$ be a finite permutation group acting transitively on the set $\Omega$ and $A_0,\dots,A_d$ be the adjacency matrices for orbits of $G$ on $\Omega\times\Omega$. 
Let $\{ B_1,\dots, B_t \} \subseteq \{ A_1,\dots, A_d \}$ be a set of adjacency matrices for the self-paired or mutually paired orbits. 
Then $\displaystyle M=\sum_{i=1}^{t}B_i$ is the adjacency matrix of a regular graph $\Gamma$. The group $G$ acts transitively on the set of vertices of the graph $\Gamma$.
\end{tm}

Using this method one can construct all regular graphs admitting a transitive action of the group $G$. We will be interested only in those 
regular graphs that are distance-regular, and especially strongly regular.

\subsection{SRGs and DRGs from the group $G_2(4)$}\label{SRG_G24}

The group $G_2(4)$ is the simple group of order $251596800=2^{12}\cdot 3^3\cdot 5^2 \cdot 7\cdot 13$. Up to conjugation it has 4300 subgroups, 8 of which are maximal. 
It belongs to Chevalley's exceptional groups of Lie type. The group $G_2(4)$ has 9 conjugacy classes of subgroups up to the index 13000. 
In Table \ref{tb:subgrpsG24} we give the list of all the subgroups $H_i^1 \le G_2(4)$ which lead to the construction of SRGs or DRGs of diameter $d \ge 3$.

\begin{table}[H]
\begin{center} \begin{footnotesize}
\begin{tabular}{|c|c|r|r|r|c|}
\hline
Subgroup& Structure& Order & Index & Rank & Primitive\\
\hline
\hline
$H^1_{1}$ & $J_2$ & 604800 & 416 & 3 & yes \\
$H^1_{2}$ & $2^{2+8}:(3\times A_5)$& 184320 & 1365 & 4 & yes \\
$H^1_{3}$ & $2^{4+6}:(A_5\times 3)$& 184320 & 1365 & 4 & yes \\
$H^1_{4}$ & $U(3,4):Z_2$&124800 & 2016 & 3 & yes \\
$H^1_{5}$ & $3^.L_3(4):2_3$& 120960 & 2080 & 4 & yes \\
$H^1_{6}$ & $2^6:(2^4:A_5)$ & 61440 & 4095 & 8 & no \\
$H^1_{7}$ & $(2^2\times((2^4:2):2)):(A_4\times A_4)$ & 36864 & 6825 & 12 & no \\
\hline
\hline
\end{tabular}\end{footnotesize} 
\caption{\footnotesize Subgroups of the group $G_2(4)$}\label{tb:subgrpsG24}
\end{center}
\end{table}

Using the method described in Theorem \ref{main} we obtained all DRGs with at most 20000 vertices on which the group $G_2(4)$ acts primitively, 
and all DRGs with at most 13000 vertices on which the group $G_2(4)$ acts imprimitively (transitively but not primitively).

Using the computer search we obtained SRGs on 416, 1365, 2016, 2080 or 4095 vertices. 
We determined the full automorphism groups of the constructed SRGs.

\begin{tm} \label{srg-G24}
Up to isomorphism there are exactly five strongly regular graphs with at most $20000$ vertices, admitting a primitive action of the group $G_2(4)$ and exactly one strongly regular graph with at most 
$13000$ vertices, admitting an imprimitive action of the group $G_2(4)$.
These strongly regular graphs have parameters $(416,100,36,20)$, $(1365,340,83,85)$, $(2080,1008,480,496)$, $(2016,975,462,480)$ and $(4095,2046,1021,1023)$. 
Details about the obtained strongly regular graphs are given in Table \ref{tb:srgG24}. 
\end{tm}

\begin{table}[H]
\begin{center} \begin{footnotesize}
\begin{tabular}{|c|c|c|}
\hline
Graph $\Gamma$ & Parameters  & $Aut (\Gamma) $  \\
\hline
\hline
$\Gamma^{1}_{1}=\Gamma(G_2(4),H^1_{1})$ & (416,100,36,20) & $G_2(4):Z_2$  \\ 
$\Gamma^{1}_{2}=\Gamma(G_2(4),H^1_{2})$ & (1365,340,83,85) & $G_2(4):Z_2$  \\ 
$\Gamma^{1}_{3}=\Gamma(G_2(4),H^1_{3})$ & (1365,340,83,85) & $O(7,4):Z_2$  \\ 
$\Gamma^{1}_{4}=\Gamma(G_2(4),H^1_{4})$ & (2016,975,462,480) & $O(7,4):Z_2$  \\ 
$\Gamma^{1}_{5}=\Gamma(G_2(4),H^1_{5})$ & (2080,1008,480,496) & $O(7,4):Z_2$  \\ 
$\Gamma^{1}_{6}=\Gamma(G_2(4),H^1_{6})$ & (4095,2046,1021,1023) &  $O(13,2)$ \\ 
\hline
\hline
\end{tabular}\end{footnotesize} 
\caption{\footnotesize SRGs constructed from the group $G_2(4)$}\label{tb:srgG24}
\end{center}
\end{table}

\begin{rem} \label{IsoSrgsG24}
The SRGs $\Gamma_2^{1}$ and $\Gamma_3^{1}$ are connected with the generalised hexagon $GH(4,4)$ of order $(4,4)$ (or shortly $H(4)$).
It is well known that the distance $3$ graph of the point graph of a generalised hexagon $GH(4,4)$ is strongly regular. 
Moreover, according to \cite[Corollary 3.5.7]{HMan} the generalised hexagon $H(q)$ is self-dual if and only if $q=3^h$. 
Therefore, the generalised hexagon $H(4)$ is not self-dual.
The graph $\Gamma_3^{1}$ is the distance $3$ graph of the point graph of the generalised hexagon $H(4)$, and 
the strongly regular graph $\Gamma_2^{1}$ is the distance $3$ graph of the point graph of the dual of the generalised hexagon $H(4)$. 
For more information see \cite{DempKantor}.
The graphs $\Gamma^{1}_{3}$ and $\Gamma^{1}_{6}$ are $O(7,4)$ and $O(13,2)$ graphs, respectively (see \cite{aeb}). 
The graphs $\Gamma^1_{1}$, $\Gamma^1_{3}$, $\Gamma^1_{4}$, $\Gamma^1_{5}$ and $\Gamma^1_{6}$ are rank $3$ graphs.
Strongly regular graph $\Gamma_1^1$ is known as $G_2(4)$ graph. It is locally the Janko graph (see \cite{pase}), and the second subconstituent of the Suzuki tower.
\end{rem}

\begin{rem} \label{desG24}
Strongly regular graphs $\Gamma^{1}_{2}$, $\Gamma^{1}_{3}$ and $\Gamma^{1}_{6}$ can be constructed from symmetric incidence matrices with all-one diagonal of symmetric block designs with parameters 
$(1365,341,85)$ and $(4095,2047,1023)$. 
Note that the full automorphism groups of the designs and the corresponding SRGs do not have to be the same. In Table \ref{tb:des} we give the details about obtained designs.

\begin{table}[H]
\begin{center} \begin{footnotesize}
\begin{tabular}{|c|c|c|}
\hline
Design $\mathcal{D}$ & Parameters  & $Aut (\mathcal{D}) $  \\
\hline
\hline
$\mathcal{D}^{1}_{1}=\Gamma(G_2(4),H^1_{2})$ & 2-(1365,341,58) & $G_2(4):Z_2$\\ 
$\mathcal{D}^{1}_{2}=\Gamma(G_2(4),H^1_{3})$ & 2-(1365,341,58) & $PSL(6,4):S_3$\\ 
$\mathcal{D}^{1}_{3}=\Gamma(G_2(4),H^1_{6})$ & 2-(4095,2047,1023) & $PSL(12,2)$ \\ \hline
\hline
\end{tabular} \end{footnotesize}
\caption{\footnotesize Symmetric block designs constructed from the group $G_2(4)$}\label{tb:des}
\end{center}
\end{table}

\end{rem}

Using the computer search we obtained distance-regular graphs on 1365 and 6825 vertices with diameter $d \geq 3$. 
We determined the full automorphism groups of the constructed DRGs. 

\begin{tm} \label{drgs-G24}
Up to isomorphism there are exactly two distance-regular graphs of diameter $d \ge 3$ with at most $20000$ vertices, 
admitting a primitive action of the group $G_2(4)$ and exactly one distance-regular graph of diameter $d \ge 3$ with at most 
$13000$ vertices, admitting an imprimitive action of the group $G_2(4)$.
These DRGs have $1365$ and $6825$ vertices, respectively. Details about the obtained DRGs are given in Table \ref{tb:drgsG24}. 
\end{tm}

\begin{table}[H]
\begin{center} \begin{footnotesize}
\begin{tabular}{|c|c|c|c|}
\hline
Graph $\Gamma$ & Number of vertices  & Intersection array & $Aut (\Gamma) $  \\
\hline
\hline
$\Gamma^{1}_{7}=\Gamma(G_2(4),H^1_{2})$ & 1365 & $\{ 20,16,16;1,1,5 \}$& $G_2(4):Z_2$\\ 
$\Gamma^{1}_{8}=\Gamma(G_2(4),H^1_{3})$ & 1365 & $\{ 20,16,16;1,1,5 \}$& $G_2(4):Z_2$\\ 
$\Gamma^{1}_{9}=\Gamma(G_2(4),H^1_{7})$ & 6825 &$\{ 8,4,4,4,4,4;1,1,1,1,1,2 \}$ & $G_2(4):Z_2$\\
\hline
\hline
\end{tabular} \end{footnotesize}
\caption{\footnotesize DRGs constructed from the group $G_2(4)$, $d \geq$ 3}\label{tb:drgsG24}
\end{center}
\end{table}

\begin{rem}
Distance-regular graphs $\Gamma^{1}_{7}$ and $\Gamma^{1}_{8}$ have diameter $3$ and $\Gamma^{1}_{9}$ has diameter $6$. The graphs $\Gamma^{1}_{7}$ and $\Gamma^{1}_{8}$ belong to the family of graphs of Lie type from Chevalley groups. 
See \cite{drg} for more information. The graph $\Gamma_9^1$ is the generalised dodecagon of order $(4,1)$. See \cite{damHaemers} for more information.
\end{rem}

\subsection{SRGs and DRGs from the group $G_2(5)$}\label{SRG_G25}

The group $G_2(5)$ is the simple group of order $ 5859000000 =2^{6}\cdot 3^5\cdot 5^6 \cdot 7\cdot 31$. Up to conjugation it has 824 subgroups, 7 of which are maximal. It belongs to Chevalley's exceptional groups of Lie type. The group $G_2(5)$ has 7 conjugacy classes of subgroups up to the index 15000. 
In Table \ref{tb:subgrpsG25} we give the list of all the subgroups $H_i^2 \le G_2(5)$ which lead to the construction of SRGs or DRGs of diameter $d \ge 3$.

\begin{table}[H]
\begin{center} \begin{footnotesize}
\begin{tabular}{|c|c|r|r|r|c|}
\hline
Subgroup& Structure& Order & Index & Rank & Primitive\\
\hline
\hline
$H^2_{1}$ & $5^{1+4}:GL_2(5)$& 1500000 & 3906 & 4 & yes \\
$H^2_{2}$ & $5^{2+3}:GL_2(5)$& 1500000 & 3906 & 4 & yes \\
$H^2_{3}$ & $3.U_3(5):2$& 756000 & 7750 & 4 & yes \\
$H^2_{4}$ & $L_3(5):2$ & 744000 & 7875 & 5 & yes \\
\hline
\hline
\end{tabular}\end{footnotesize} 
\caption{\footnotesize Subgroups of the group $G_2(5)$}\label{tb:subgrpsG25}
\end{center}
\end{table}

Using the method described in Theorem \ref{main} we obtained all DRGs with at most 400000 vertices on which the group $G_2(5)$ acts primitively, and all DRGs with at most 15000 vertices 
on which the group $G_2(5)$ acts imprimitively.

Using the computer search we obtained SRGs on 3906, 7750 or 7875 vertices. 
Further, we determined the full automorphism groups of the constructed SRGs.

\begin{tm} \label{Srg-G25}
Up to isomorphism there are exactly $6$ strongly regular graphs with at most $400000$ vertices, admitting a primitive action of the group $G_2(5)$.
These strongly regular graphs have parameters $(3906,780,154,156)$, $(7875,1550,325,300)$, $(7875,3224,1348,1300)$, $(7750,3024,1148,1200)$ and $(7750,1575,300,325)$. 
Details about the obtained strongly regular graphs are given in Table \ref{tb:srgG25}. 
\end{tm}


\begin{table}[H]
\begin{center} \begin{footnotesize}
\begin{tabular}{|c|c|c|}
\hline
Graph $\Gamma$ & Parameters  & $Aut (\Gamma) $  \\
\hline
\hline
$\Gamma^{2}_{1}=\Gamma(G_2(5),H^2_{1})$ & $(3906,780,154,156)$ & $G_2(5)$  \\ 
$\Gamma^{2}_{2}=\Gamma(G_2(5),H^2_{2})$ & $(3906,780,154,156)$ & $O(7,5):Z_2$  \\ 
$\Gamma^{2}_{3}=\Gamma(G_2(5),H^2_{3})$ & $(7750,1575,300,325)$ & $O(7,5)$  \\ 
$\Gamma^{2}_{4}=\Gamma(G_2(5),H^2_{3})$ & $(7750,3024,1148,1200)$ & $O(7,5)$  \\ 
$\Gamma^{2}_{5}=\Gamma(G_2(5),H^2_{4})$ & $(7875,3224,1348,1300)$ & $O(7,5)$  \\ 
$\Gamma^{2}_{6}=\Gamma(G_2(5),H^2_{4})$ & $(7875,1550,325,300)$ & $O(7,5)$  \\ 
\hline
\hline
\end{tabular}\end{footnotesize} 
\caption{\footnotesize SRGs constructed from the group $G_2(5)$}\label{tb:srgG25}
\end{center}
\end{table}

\begin{rem} \label{IsoSrgsG25}
The SRGs $\Gamma_1^{2}$ and $\Gamma_2^{2}$ are connected with the generalised hexagon $GH(5,5)$ of order $(5,5)$ (or shortly $H(5)$).
The graph $\Gamma_2^{2}$ is the distance $3$ graph of the point graph of the generalised hexagon $H(5)$, and 
the strongly regular graph $\Gamma_1^{2}$ is the distance $3$ graph of the point graph of the dual of the generalised hexagon $H(5)$.
The graph $\Gamma^{2}_{2}$ is a rank $3$ graph and known as the $O(7,5)$ graph.
Further, strongly regular graphs $\Gamma^{2}_{1}$ and $\Gamma^{2}_{2}$ can be constructed from a symmetric incidence matrix with all-one diagonal of a symmetric block design with parameters 
$2$-$(3906,781,156)$. 
The full automorphism groups of symmetric block designs in these cases are the same as the automorphism groups of the corresponding SRGs.
The SRGs with parameters $(7875,1550,325,300)$, $(7875,3224,1348,1300)$, $(7750,3024,1148,1200)$ and $(7750,1575,300,325)$ belong to families of SRGs related to non-singular quadrics in $PG(2m,q)$, 
described in \cite[Chapters 7C,7D]{LintBrouwer}, for $m=3$ and $q=5$.
We checked imprimitive representations of the group $G_2(5)$ up to the degree $15000$, and proved that there are no SRGs in these cases.
\end{rem}

Using the computer search we obtained a distance-regular graph on 3906 vertices with diameter $d \geq 3$. 
We determined the full automorphism group of the constructed DRG. 

\begin{tm} \label{drgs-G25}
Up to isomorphism there are exactly two distance regular graphs of diameter $d \ge 3$ with at most $400000$ vertices, admitting a primitive action of the group $G_2(5)$.
They have $3906$ vertices. Details about the obtained DRGs are given in Table \ref{tb:drgG25}. 
\end{tm}

\begin{table}[H]
\begin{center} \begin{footnotesize}
\begin{tabular}{|c|c|c|c|}
\hline
Graph $\Gamma$ & Number of vertices  & Intersection array & $Aut (\Gamma) $  \\
\hline
\hline
$\Gamma^{2}_{7}=\Gamma(G_2(5),H^2_{1})$ & 3906 & $\{ 30,25,25; 1,1,6 \}$& $G_2(5)$\\ 
$\Gamma^{2}_{8}=\Gamma(G_2(5),H^2_{2})$ & 3906 & $\{ 30,25,25; 1,1,6 \}$& $G_2(5)$\\ 
\hline
\hline
\end{tabular} \end{footnotesize}
\caption{\footnotesize DRG constructed from the group $G_2(5)$, $d \geq$ 3}\label{tb:drgG25}
\end{center}
\end{table}

\begin{rem}
The distance-regular graphs $\Gamma^{2}_{7}$ and $\Gamma^{2}_{8}$ have diameter $3$. They belong to the family of graphs of Lie type from Chevalley groups. See \cite{drg} for more information.
We checked imprimitive representations of the group $G_2(5)$ up to the degree $15000$, and proved that there are no DRGs in these cases.
\end{rem}

\subsection{SRGs and DRGs from the group $O(7,3)$}\label{SrgsO73}

The group $O(7,3)$ is the simple group of order $4585351680=2^{9}\cdot 3^9\cdot 5 \cdot 7\cdot 13$. Up to conjugation it has 7735 subgroups, 15 of which are maximal. The group $O(7,3)$ has 15 conjugacy classes of subgroups up to the index 10000. 
In Table \ref{tb:subgrpsO73} we give the list of all the subgroups $H_i^3 \le O(7,3)$ which lead to the construction of SRGs or DRGs of diameter $d \ge 3$.

\begin{table}[H]
\begin{center} \begin{footnotesize}
\begin{tabular}{|c|c|r|r|r|c|}
\hline
Subgroup& Structure& Order & Index & Rank & Primitive\\
\hline
\hline
$H^3_{1}$ & $2U_4(3):2$ & 13063680 & 351 & 3 & yes \\
$H^3_{2}$ & $3^5:U_4(2):2$& 12597120 & 364 & 3 & yes \\
$H^3_{3}$ & $L_4(3):2$& 12130560 & 378 & 3 & yes \\
$H^3_{4}$ & $G_2(3)$&4245696 & 1080 & 3 & yes \\
$H^3_{5}$ & $3^{3+3}:L_3(3)$& 4094064 & 1120 & 4 & yes \\
$H^3_{6}$ & $2^6:A_7$ & 161280 & 28431 & 12 & yes \\
\hline
\hline
\end{tabular}\end{footnotesize} 
\caption{\footnotesize Subgroups of the group $O(7,3)$}\label{tb:subgrpsO73}
\end{center}
\end{table}

Using the method described in Theorem \ref{main} we obtained all DRGs with at most 250000 vertices on which the group $O(7,3)$ acts primitively, and all DRGs with at most 10000 vertices 
on which the group $O(7,3)$ acts imprimitively.

Using the computer search we obtained SRGs on 351, 364, 378, 1080, 1120 or 28431 vertices. 
We determined the full automorphism groups of the constructed SRGs.

\begin{tm} \label{srg-O73}
Up to isomorphism there are exactly $7$ SRGs with at most $250000$ vertices, admitting a primitive action of the orthogonal group $O(7,3)$.
These SRGs have parameters $(351,126,45,45)$, $(364,120,38,40)$, $(378,117,36,36)$, $(1080,351,126,108)$, $(1120,390,146,130)$, $(28431,2880,324,288)$ and $(28431,3150,621,315)$. 
Details about the obtained SRGs are given in Table \ref{tb:srgO73}. 
\end{tm}


\begin{table}[H]
\begin{center} \begin{footnotesize}
\begin{tabular}{|c|c|c|}
\hline
Graph $\Gamma$ & Parameters  & $Aut (\Gamma) $  \\
\hline
\hline
$\Gamma^{3}_{1}=\Gamma(O(7,3),H^3_{1})$ & (351,126,45,45) & $O(7,3):Z_2$  \\ 
$\Gamma^{3}_{2}=\Gamma(O(7,3),H^3_{2})$ & (364,120,38,40) & $O(7,3):Z_2$  \\ 
$\Gamma^{3}_{3}=\Gamma(O(7,3),H^3_{3})$ & (378,117,36,36) & $O(7,3):Z_2$  \\ 
$\Gamma^{3}_{4}=\Gamma(O(7,3),H^3_{4})$ & (1080,351,126,108) & $O^{+}(8,3):E_4$  \\ 
$\Gamma^{3}_{5}=\Gamma(O(7,3),H^3_{5})$ & (1120,390,146,130) & $O^{+}(8,3):D_8$  \\ 
$\Gamma^{3}_{6}=\Gamma(O(7,3),H^3_{6})$ & (28431,3150,621,315) &  $O^{+}(8,3):S_3$  \\ 
$\Gamma^{3}_{7}=\Gamma(O(7,3),H^3_{6})$ & (28431,2880,324,288) &  $O^{+}(8,3):S_3$  \\ 
\hline
\hline
\end{tabular}\end{footnotesize} 
\caption{\footnotesize SRGs constructed from the group $O(7,3)$}\label{tb:srgO73}
\end{center}
\end{table}

\begin{rem} \label{IsoSrgsO73}
To the best of our knowledge the SRGs with parameters $(28431,3150,621,315)$ and $(28431,2880,324,288)$ are the first known examples of SRGs with these parameters. 
The graphs $\Gamma^3_{1}$, $\Gamma^3_{2}$, $\Gamma^3_{3}$, $\Gamma^3_{4}$ and $\Gamma^3_{5}$ are rank $3$ graphs. The graph $\Gamma^3_{1}$ is $NO^{-1}(7,3)$ graph, $\Gamma^3_{2}$ is $O(7,3)$ graph, 
the graph $\Gamma^3_{3}$ is $NO^{+1}(7,3)$ graph, the graph $\Gamma^3_{4}$ is $NO^{+}(8,3)$ graph and the graph $\Gamma^3_{5}$ is $O^{+}(8,3)$ graph.
The orthogonal group $O^+(8,3)$ acts primitively on the SRGs $\Gamma^3_{4}$, $\Gamma^3_{5}$, $\Gamma^3_{6}$ and $\Gamma^3_{7}$.
Further, strongly regular
graphs $\Gamma^{3}_{1}$ and $\Gamma^{3}_{3}$ can be constructed from symmetric incidence matrices with all-zero diagonal of symmetric block designs with parameters $2$-$(351,126,45)$ and 
$2$-$(378,117,36)$ respectively. Similarly, strongly regular
graph $\Gamma^{3}_{2}$ can be constructed from symmetric incidence matrix with all-one diagonal of a symmetric block design with parameters $2$-$(364,121,40)$. 
The full automorphism groups of symmetric block designs in these cases are the same as the automorphism groups of the corresponding SRGs.
We checked imprimitive representations of the group $O(7,3)$ up to the degree $10000$, and proved that there are no SRGs in these cases.
\end{rem}

Using the computer search we obtained distance-regular graph on 1120 vertices with diameter $d \geq 3$. 
Finally, we determined the full automorphism group of the constructed DRG. 

\begin{tm} \label{drgsO73}
Up to isomorphism there is exactly one distance regular graph of diameter $d \ge 3$ with at most $250000$ vertices, admitting a primitive action of the group $O(7,3)$.
It has $1120$ vertices. Details about the obtained DRG are given in Table \ref{tb:drgO73}. 
\end{tm}

\begin{table}[H]
\begin{center} \begin{footnotesize}
\begin{tabular}{|c|c|c|c|}
\hline
Graph $\Gamma$ & Number of vertices  & Intersection array & $Aut (\Gamma) $  \\
\hline
\hline
$\Gamma^{3}_{8}=\Gamma(O(7,3),H^3_{5})$ & 1120 & $\{ 39,36,27; 1,4,13 \}$& $O(7,3):Z_2$\\ 
\hline
\hline
\end{tabular} \end{footnotesize}
\caption{\footnotesize DRG constructed from the group $O(7,3)$, $d \geq$ 3}\label{tb:drgO73}
\end{center}
\end{table}

\begin{rem}
The graph $\Gamma_8^3$ is the dual polar graph $B_3(3)$ (see \cite{drg}).
We checked imprimitive representations of the group $O(7,3)$ up to the degree $10000$, and proved that there are no DRGs in these cases.
\end{rem}

\subsection{SRGs and DRGs from the Tits group $T$}\label{SRG_T}

The Tits group $^2F_4(2)'$ or $T$, named after Jacques Tits, is a finite simple group of order $17971200=2^{11}\cdot 3^3\cdot 5^2\cdot 13$. Up to conjugation it has 434 subgroups, 8 of which are maximal.
It is sometimes considered as the $27^{\mathrm{th}}$ sporadic group since the Tits group itself is not a group of Lie type even the group $^2F_4(2)$ (which is not simple group) is a group of Lie type.
The Tits group $T$ has 10 conjugacy classes of subgroups up to the index 7000. 
In Table \ref{tb:subgrpsT} we give the list of all the subgroups $H_i^4 \le T$ which lead to the construction of SRGs or DRGs of diameter $d \ge 3$.

\begin{table}[H]
\begin{center} \begin{footnotesize}
\begin{tabular}{|c|c|r|r|r|c|}
\hline
Subgroup& Structure& Order & Index & Rank & Primitive\\
\hline
\hline
$H^4_{1}$ & $L_3(3):2$& 11232 & 1600 & 4 & yes \\
$H^4_{2}$ & $2.2^8:5:4$& 10240 & 1755 & 5 & yes \\
$H^4_{3}$ & $2.2^8:S_3$& 6144 & 2925 & 6 & yes \\
$H^4_{4}$ & $2^5:((2^4:5):2)$& 5120 & 3510 & 10 & no \\
\hline
\hline
\end{tabular}\end{footnotesize} 
\caption{\footnotesize Subgroups of the group $T$}\label{tb:subgrpsT}
\end{center}
\end{table}

Using the method described in Theorem \ref{main} we obtained all DRGs with at most 12000 vertices on which the group $T$ acts primitively, and all DRGs with at most 7000 vertices on which the group 
$T$ acts imprimitively.

Using the computer search we obtained SRGs on 1600 or 3510 vertices. 
Finally, we determined the full automorphism groups of the constructed SRGs.

\begin{tm} \label{srg-G25}
Up to isomorphism there is exactly one strongly regular graph with at most $12000$ vertices, admitting a primitive action of the group $T$ and exactly one strongly regular graph with at most $7000$ vertices, admitting an imprimitive action of the group $T$.
These strongly regular graphs have parameters  $(1600,351,94,72)$ and $(3510,693,180,126)$. Details about the obtained strongly regular graphs are given in Table \ref{tb:srgT}. 
\end{tm}

\begin{table}[H]
\begin{center} \begin{footnotesize}
\begin{tabular}{|c|c|c|}
\hline
Graph $\Gamma$ & Parameters  & $Aut (\Gamma) $  \\
\hline
\hline
$\Gamma^{4}_{1}=\Gamma(T,H^4_{1})$ & $(1600,351,94,72)$ & $T$  \\ 
$\Gamma^{4}_{2}=\Gamma(T,H^4_{4})$ & $(3510,693,180,126)$ & $Fi_{22}.2$  \\ 

\hline
\hline
\end{tabular}\end{footnotesize} 
\caption{\footnotesize SRGs constructed from the group $T$}\label{tb:srgT}
\end{center}
\end{table}

\begin{rem} \label{IsoSrgsT}
The SRGs with parameters $(1600,351,94,72)$ and $(3510,693,180,126)$ are described in \cite{ieee}.
The graph $\Gamma^{4}_{2}$ can be constructed from the group $Fi_{22}$ as a rank $3$ graph. Moreover, it belongs to the Fisher tower (see \cite{hubaut, tits}).
\end{rem}

Using the computer search we obtained DRGs on 1755 or 2925 vertices with diameter $d \geq 3$. 
Finally, we determined the full automorphism group of the constructed DRGs. 

\begin{tm} \label{drgsT}
Up to isomorphism there are exactly two distance-regular graphs of diameter $d \ge 3$ with at most $12000$ vertices, admitting a primitive action of the group $T$.
These DRGs have $1755$ or $2925$ vertices. Details about the obtained DRGs are given in Table \ref{tb:drgT}. 
\end{tm}

\begin{table}[H]
\begin{center} \begin{footnotesize}
\begin{tabular}{|c|c|c|c|}
\hline
Graph $\Gamma$ & Number of vertices  & Intersection array & $Aut (\Gamma) $  \\
\hline
\hline
$\Gamma^{4}_{3}=\Gamma(T,H^4_{2})$ & 1755 & $\{ 10,8,8,8; 1,1,1,5 \}$& $T:2$\\ 
$\Gamma^{4}_{4}=\Gamma(T,H^4_{3})$ & 2925 & $\{ 12,8,8,8; 1,1,1,3 \}$& $T:2$\\ 

\hline
\hline
\end{tabular} \end{footnotesize}
\caption{\footnotesize DRGs constructed from the group $T$, $d \geq$ 3}\label{tb:drgT}
\end{center}
\end{table}

\begin{rem}
Distance-regular graphs $\Gamma^{4}_{3}$ and $\Gamma^{4}_{4}$ have diameter $4$. The graph $\Gamma^{4}_{3}$ is the generalized octagon of order $(2,4)$ and the graph 
$\Gamma^{4}_{4}$ is the generalized octagon of order $(4,2)$.
We checked imprimitive representations of the group $T$ up to the degree $7000$, and proved that there are no DRGs in these cases.
\end{rem}

\begin{rem}
All SRGs and DRGs constructed in this paper are vertex-transitive, but some of them are also edge-transitive. 
The graphs that are not edge-transitive are $\Gamma^1_{2}$, $\Gamma^1_{8}$ and $\Gamma^2_{1}$.
\end{rem}

\vspace*{0.2cm}

\noindent {\bf Acknowledgement} \\
This work has been fully supported by {\rm C}roatian Science Foundation under the project 6732. 
The authors would like to thank Ferdinand Ihringer, William M. Kantor and Dmitrii Pasechnik for their valuable suggestions.

\end{document}